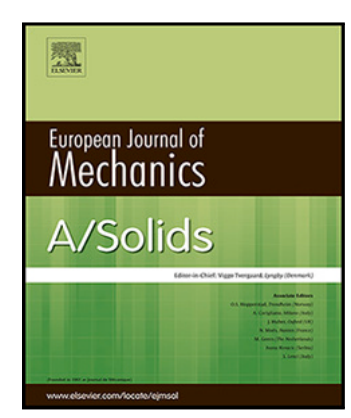

# Polyconvexity implies Hill's inequality in SL(2)

Ionel-Dumitrel Ghiba [a,b,*], Maximilian P. Wollner [c], Patrizio Neff [d]

[a] Department of Mathematics, Alexandru Ioan Cuza University of Iaşi, Blvd. Carol I, no. 11, 700506 Iaşi, Romania
[b] Octav Mayer Institute of Mathematics of the Romanian Academy, Iaşi Branch, 700505 Iaşi, Romania
[c] Institute of Biomechanics, Graz University of Technology, Stremayrgasse 16/2, 8010, Graz, Austria
[d] Head of Lehrstuhl für Nichtlineare Analysis und Modellierung, Fakultät für Mathematik, Universität Duisburg–Essen, Thea-Leymann Str. 9, 45127 Essen, Germany

ABSTRACT

For compressible nonlinear isotropic elasticity it is well known that rank-one convexity, polyconvexity and the monotonicity of the Cauchy stress tensor with respect to the logarithmic stretch tensor (the true-stress–true-strain monotonicity, TSTS-M$^+$) are independent constitutive conditions which should, however, all together be satisfied for a physically meaningful description of idealized elastic materials. In the incompressible case, TSTS-M$^+$ turns into Hill's inequality since the Cauchy stress $\boldsymbol{\sigma}$ reduces to the Kirchhoff stress $\boldsymbol{\tau}$. Hill's inequality requires then monotonicity of the Kirchhoff stress in terms of the logarithmic stretch tensor evaluated for incompressible response. In this paper we clarify how the a priori independent notions of Legendre–Hadamard ellipticity (LH), polyconvexity and Hill's inequality are nevertheless intimately connected. More precisely, by providing several alternative proofs, we show that both LH-ellipticity (rank-one convexity) and polyconvexity imply the weak Hill inequality in the incompressible two-dimensional case.

## 1. Introduction

Incompressible nonlinear elasticity is perhaps the most widely used assumption in the analytical and algorithmic treatment of elastomeric materials. While these materials are indeed highly compressible, they still show a comparatively large bulk modulus (high resistance against volumetric deformation) compared to a small shear modulus (low resistance against shape change), so that for certain boundary value problems, the actual response is in fact nearly incompressible (isochoric) (Treloar, 1975). This incompressibility assumption can be used to good advantage in deriving analytical solutions to certain non-trivial boundary value problems, since PDE-systems can then be reduced to ODE-questions. On the other hand, material stability questions need to be discussed in advance before setting up elastic energies. Among them, the Legendre–Hadamard ellipticity condition (Knowles and Sternberg, 1975; Dacorogna, 2008) and Ball's polyconvexity condition (Ball, 1977) have emerged as important concepts.

Recently, it has been shown that polyconvexity alone is, however, insufficient to guarantee, e.g., monotone increasing Cauchy stress in uniaxial tension with free lateral sides (Wollner et al., 2026a). Insisting on a multiaxial version of Cauchy stress increases with strain, Neff and coauthors (Neff et al., 2025a,b) have put forward the true-stress–true-strain monotonicity condition (TSTS-M$^+$), which demands in terms of the Hencky strain tensor $\log V$ (Jog and Patil, 2013; Neff et al., 2016, 2015a):

$$\langle \boldsymbol{\sigma}(\log V_1) - \boldsymbol{\sigma}(\log V_2), \log V_1 - \log V_2 \rangle > 0 \quad \forall\, V_1, V_2 \in \mathrm{Sym}^{++}(n), \tag{TSTS-M$^+$}$$

and which implies, among others, increasing Cauchy-stress in uniaxial tension. It can be shown that the TSTS-M$^+$ condition is equivalent (Martin and Neff, 2026) to the so called corotational stability postulate (CSP), meaning that

$$\left\langle \frac{\mathrm{D}^\circ}{\mathrm{D}t}[\boldsymbol{\sigma}],\, D \right\rangle > 0 \qquad \text{for all} \quad D = \mathrm{sym}(L), \tag{CSP}$$

where $L = \dot{F}F^{-1}$ is the velocity gradient (Neff et al., 2025a) and $\frac{\mathrm{D}^\circ}{\mathrm{D}t}$ is any suitable corotational stress rate (Leblond, 1992; Neff et al., 2026, 2024).

Furthermore, it is known that TSTS-M$^+$ and polyconvexity are essentially independent constitutive requirements (Wollner et al., 2026a; Leblond, 1992) as neither implies the other.

For incompressibility, however, the situation changes considerably. First, the TSTS-M$^+$ condition turns into the better known Hill's inequality, inasmuch as the Cauchy stress $\boldsymbol{\sigma} = \frac{1}{\det F}\boldsymbol{\tau}$ reduces to the Kirchhoff stress $\boldsymbol{\tau}$, see e.g., (Hill, 1957, 1968b,a, 1970; Sidoroff, 1974)

$$\langle \boldsymbol{\tau}(V_1) - \boldsymbol{\tau}(V_2), \log V_1 - \log V_2 \rangle > 0 \quad \forall\, V_1, V_2 \in \mathrm{Sym}^{++}(n). \tag{Hill}$$

* Corresponding author at: Department of Mathematics, Alexandru Ioan Cuza University of Iaşi, Blvd. Carol I, no. 11, 700506 Iaşi, Romania.
*E-mail addresses:* dumitrel.ghiba@uaic.ro (I.-D. Ghiba), wollner@tugraz.at (M.P. Wollner), patrizio.neff@uni-due.de (P. Neff).

Due to incompressibility, we have $\det F = 1 \iff \det V = 1 \iff \operatorname{tr}(\log V) = 0$, where $V = \sqrt{F F^T}$ is the left stretch tensor, so that the monotonicity condition (Hill) needs only to be evaluated over the linear constraint set $\operatorname{tr}(\log V) = 0$. We will refer to this inequality loosely as the ***weak Hill inequality***.

In the isotropic case, for smooth isotropic and objective functions $W : \mathrm{GL}^+(n) \to \mathbb{R}$ we have

$$W(F) = g(\lambda_1, \lambda_2, \ldots, \lambda_n) = \widehat{g}(\log\lambda_1, \log\lambda_2, \ldots, \log\lambda_n) = \widehat{W}(\log V), \qquad (1.1)$$

where $\lambda_1, \lambda_2, \ldots, \lambda_n$ are the singular values of $F$, $g : \mathbb{R}^n_+ \to \mathbb{R}$ and $\widehat{g} : \mathbb{R}^n \to \mathbb{R}$ are permutation-symmetric and $\widehat{W}(\log V) : \mathrm{Sym}(n) \to \mathbb{R}$. Moreover, it holds (see, e.g., Richter (1949, 1948), Neff et al. (2020)) that the Kirchhoff stress tensor $\boldsymbol{\tau}$ has an elastic potential

$$\boldsymbol{\tau} = \mathrm{D}_{\log V}\, \widehat{W}(\log V) \qquad (1.2)$$

in terms of the logarithmic strain $\log V$. Therefore, the Hill monotonicity condition translates into a convexity requirement for the mapping

$$\log V \;\mapsto\; \widehat{W}(\log V) \qquad (1.3)$$

on the set

$$\mathrm{SL}(n) \cap \mathrm{Sym}^{++}(n) := \{\, V \in \mathrm{Sym}^{++}(n) \mid \det V = 1 \,\}. \qquad (1.4)$$

Thus, we have that $\widehat{W}(X)$ needs to be convex only over the linear constraint set $\operatorname{tr}(X) = 0$. In Baaser (2026, Appendix-D3-Neff) the appropriate reduced convexity condition has been first considered for $n = 3$. This is most easily achieved by switching to the representation in principal logarithmic stretches. Indeed, for smooth functions $W(F) = g^{\rm inc}(\lambda_1, \lambda_2, \lambda_3) = \widehat{g}^{\rm inc}(\log\lambda_1, \log\lambda_2, \log\lambda_3)$, we may define

$$\begin{aligned} W^{\rm inc}_{\rm red}(\lambda_1, \lambda_2) &:= g^{\rm inc}\left(\lambda_1, \lambda_2, \frac{1}{\lambda_1\lambda_2}\right) := \widehat{g}^{\rm inc}_{\rm red}(\log\lambda_1, \log\lambda_2) \\ &= \widehat{g}^{\rm inc}\big(\log\lambda_1, \log\lambda_2, -(\log\lambda_1 + \log\lambda_2)\big)\,. \end{aligned} \qquad (1.5)$$

and the weak Hill inequality is satisfied if and only if $\widehat{g}^{\rm inc}_{\rm red}$ is convex.

In Wollner et al. (2026b) we investigated the concurrent enforcement of polyconvexity and TSTS–M$^+$ for the incompressible case, together with its implementation to Physics Informed Neural Networks (PINN). It turned out that a wide range of polyconvex functions already satisfy the weak Hill inequality, e.g., all polyconvex energies which satisfy Ball's sufficient condition (Ball, 1977) (see also the paper by Steigmann (Steigmann, 2003), page 485) and moreover all incompressible energies[1] $W = W(I_1, I_2)$, convex and monotonically increasing in $I_1 = \operatorname{tr} C, I_2 = \operatorname{tr}\operatorname{Cof} C$, $C = F^T F$ are both polyconvex and satisfy the weak Hill inequality. The polyconvexity condition is often taken as the defining requirement for material stability (see, e.g., Poya et al. (2025)).

To the best of our knowledge, the direct relationship between Hill's inequality and polyconvexity has not previously been systematically investigated as an implication problem. Recent work by Wollner and coauthors (Wollner et al., 2026b,a) has considered the simultaneous enforcement of polyconvexity and true-stress–true-strain monotonicity and has identified constitutive conditions and broad classes of energy functions for which both requirements are satisfied. These results show that the two conditions are compatible for important classes of constitutive models, but they do not establish a general implication between them. The present contribution addresses precisely this question in the incompressible two-dimensional setting and shows that, due to the particular structure of objective and isotropic energies on SL(2), polyconvexity (equivalently, rank-one convexity) implies the weak Hill inequality. Thus, our result reveals a previously unnoticed connection between constitutive requirements which, in the general nonlinear elastic setting, are regarded as independent.

In this paper we reconsider the situation for the weak Hill inequality in the smooth and non-smooth case, where we first revisit some results related to polyconvexity in SL(2).

In Section 2, we give an overview of the main constitutive requirements which we use in the present paper, namely: rank-one convexity, polyconvexity, and Hill's inequality, as well as a short overview of some theorems characterizing these constitutive requirements. Let us mention that all these constitutive properties are considered in this paper for energies defined on SL($n$); we consider their intrinsic definitions on this manifold. We give special attention to the definition of the Kirchhoff stress tensor $\boldsymbol{\tau}$ for functions on SL($n$) and we show the relation between its matrix monotonicity in $\log V$ and the vector monotonicity in terms of its principal components $\log\lambda_i$.

Then we prove the main result of the present paper:

**Theorem 1.2** (*Polyconvexity Implies Hill's Inquality in* SL(2) *for Differentiable Energies*)**.** *Let $W^{\rm inc} : \mathrm{SL}(2) \to \mathbb{R}$ be objective and isotropic, and assume that $W^{\rm inc}$ is polyconvex and differentiable. Then $W^{\rm inc}$ satisfies the weak Hill inequality.*

There are two main ingredients for the proof of this result. The first one is the result given in Ghiba et al. (2018), showing that rank-one convexity and polyconvexity for energies defined on SL(2) are equivalent. The second ingredient is the equivalence between matrix and vector monotonicity of the Kirchhoff stress tensor defined by differentiable functions on SL($n$); the latter useful result will also be proven later. In addition, the necessary and sufficient criteria obtained by Wiedemann and Peter in Wiedemann and Peter (2026) allow for an alternative proof of Theorem 1.2.

## 2. Convexity concepts and Hill's inequality on SL($n$)

### 2.1. Rank-one convexity and polyconvexity on SL($n$)

Throughout this article, we will use the following definitions of rank-one convexity, polyconvexity and matrix and vector monotonicity for energy functions defined on the sets on $\mathbb{R}^{n\times n}$, $\mathrm{GL}^+(n) := \{X \in \mathbb{R}^{n\times n} \mid \det X > 0\}$ and $\mathrm{SL}(n) := \{X \in \mathbb{R}^{n\times n} \mid \det X = 1\}$, respectively, where $n \geq 1$. The elastic energy $W^{\rm inc} : \mathrm{SL}(n) \to \mathbb{R}$ is assumed to be *objective* as well as *isotropic*, i.e., it has to satisfy the condition

$$W^{\rm inc}(Q_1\, F\, Q_2) = W^{\rm inc}(F) \quad \text{for all } F \in \mathrm{SL}(n) \text{ and all } Q_1, Q_2 \in \mathrm{SO}(n)\,, \qquad (2.1)$$

where $\mathrm{SO}(n) = \{X \in \mathbb{R}^{n\times n} \mid X^T X = \mathbb{1}\,,\ \det(X) = 1\}$ denotes the special orthogonal group. We denote by $\lambda_i$, $i = 1, 2, \ldots, n$ the singular values of $F$ (i.e. the eigenvalues of $U = \sqrt{F^T F}$), and $\lambda_{\max} := \max_{i=1,2,\ldots,n} \lambda_i$ denotes the largest singular value of $F$ (also called the *spectral norm* of $F$).

Therefore, there exists a function $g^{\rm inc} : \mathcal{S} \to \mathbb{R}$, such that

$$W^{\rm inc}(F) = g^{\rm inc}(\lambda_1, \lambda_2, \ldots, \lambda_n), \qquad (2.2)$$

where $g^{\rm inc}$ is permutation-symmetric and

$$\mathcal{S} := \{(\lambda_1, \lambda_2, \ldots, \lambda_n) \in \mathbb{R}^n_+ \mid \lambda_1\, \lambda_2\, \ldots\, \lambda_n = 1\}. \qquad (2.3)$$

[1] In terms of the invariants $(I_1, I_2, I_3)$, the following sufficient criterion for polyconvexity holds true:

**Theorem 1.1** Ball's sufficient condition (Ball, 1977) (see also the paper by Steigmann (2003), page 485) Suppose that

(i) $\psi(I_1, I_2, I_3)$ is a convex function of $(I_1, I_2, I_3)$ jointly, and

(ii) $\psi(I_1, I_2, I_3)$ is a non-decreasing function of $I_1$ and $I_2$, separately.

Then $W(F) = \psi(I_1, I_2, I_3)$ is polyconvex.

Here, $I_3 = \det(C)$ and the domain over which $\psi(I_1, I_2, I_3)$ is defined is given by $\lambda_1 > 0$, $\lambda_2 > 0$, $\lambda_3 > 0$, i.e., given by the equation $\lambda^3 - I_1\lambda^2 + I_2\lambda - I_3 = 0$ with three positive real solutions. But this domain is not convex. Therefore, it would be more adequate to say that $\psi(I_1, I_2, I_3)$ is convex in the sense of Busemann, Ewald and Shephard's definition (Busemann et al., 1963), i.e. $\psi$ can be extended to a convex function defined on the convex hull of its domain of definition.

Note that if $W^{\rm inc} : \mathrm{SL}(n) \to \mathbb{R}$ is the restriction of an isotropic and objective energy function $W : \mathrm{GL}^+(n) \to \mathbb{R}$ (or of $\mathbb{W} : \mathbb{R}^{n\times n} \to \mathbb{R}$) having the representation

$$W(F) = g(\lambda_1, \lambda_2, \dots, \lambda_n) \qquad (\text{or} \quad \mathbb{W}(F) = \mathfrak{g}(\lambda_1, \lambda_2, \dots, \lambda_n)), \tag{2.4}$$

with $g : \mathbb{R}^n_+ \to \mathbb{R}$, then $g^{\rm inc}$ is the restriction of $g$ on $\mathcal{S}$.

**Definition 2.1.** A function $W^{\rm inc} : \mathrm{SL}(n) \to \mathbb{R}$ is called *rank-one convex* if the mapping

$$t \mapsto W^{\rm inc}(F + t\,\xi \otimes \eta) \tag{2.5}$$

is convex on $\mathbb{R}$ for all $F \in \mathrm{SL}(n)$ and all $\xi, \eta \in \mathbb{R}^2$ such that $\xi \otimes \eta \in T_{\mathrm{SL}(n)}(F)$, where $T_{\mathrm{SL}(n)}(F)$ is the tangent space to $\mathrm{SL}(n)$ at $F \in \mathrm{SL}(n)$.

Since $\mathrm{rank}(\xi \otimes \eta) = 1$ implies $\det(\xi \otimes \eta) = 0$, we have[2]

$$\det(F + t\,\xi \otimes \eta) = \det F\,(1 + t\,\langle F^{-T}, \xi \otimes \eta\rangle) \tag{2.6}$$

for $F \in \mathbb{R}^{ntimesn}$. In particular, $\det(F + t\,\xi \otimes \eta) = 1$ is satisfied if and only if

$$\langle F^{-T}, \xi \otimes \eta\rangle = 0 \qquad \Longleftrightarrow \qquad \langle F^{-1}\xi, \eta\rangle = 0. \tag{2.7}$$

In addition, the general relation $\mathrm{D}_F(\det F).H = (\det F)\langle F^{-T}, H\rangle$ implies that the tangent space $T_{\mathrm{SL}(n)}$ to $\mathrm{SL}(n)$ at $F \in \mathrm{SL}(n)$ is given by

$$T_{\mathrm{SL}(n)}(F) := \{H \in \mathbb{R}^{n\times n} \mid \langle H, F^{-T}\rangle = 0\}\,. \tag{2.8}$$

**Definition 2.2** ((Ball, 1977; Mielke, 2005)). A function $W^{\rm inc} : \mathrm{SL}(2) \to \mathbb{R}$ is called polyconvex if there exists a convex function $P : \mathbb{R}^{n\times n} \times \mathbb{R} \to \mathbb{R} \cup \{\infty\}$ such that

$$W^{\rm inc}(F) = P(F, 1) \qquad \text{for all} \quad F \in \mathrm{SL}(n). \tag{2.9}$$

The rank-one convexity and polyconvexity for energies defined on SL(2) is completely characterized by the following theorem:

**Theorem 2.3** (Ghiba et al., 2018). *Let $W^{\rm inc} : \mathrm{SL}(2) \to \mathbb{R}$ be an objective and isotropic function. Then the following are* **equivalent***:*

*(i) $W^{\rm inc}$ is rank-one convex,*
*(ii) $W^{\rm inc}$ is polyconvex,*
*(iii) the function $\widetilde{\phi} : \mathbb{R} \to \mathbb{R}$ with $W^{\rm inc}(F) = \widetilde{\phi}(\lambda_1 - \lambda_2)$ for all $F \in \mathrm{SL}(2)$ with singular values $\lambda_1, \lambda_2$ is convex,*
*(iv) the function $\phi : [0, \infty) \to \mathbb{R}$ with $W^{\rm inc}(F) = \phi(\sqrt{\|F\|^2 - 2}) = \phi\Big(\lambda_{max}(F) - \dfrac{1}{\lambda_{max}(F)}\Big)$ is nondecreasing and convex, where $\lambda_{max}(F)$ denotes the maximum singular value of $F$.*

The characterization of polyconvex energies on SL(2) by criterion iv) in Theorem 2.3 is originally due to Mielke (2005). A criterion for the rank-one convexity of a twice differentiable energy in terms of $\|F\|^2 = \lambda_1^2 + \lambda_2^2$ has previously been given by Abeyaratne (1980).

Since for the polyconvexity of a function defined on SL(2) it is enough to show that its extension

$$\mathbb{W} : \mathbb{R}^{2\times 2} \to \mathbb{R} \cup \{\infty\}\,, \qquad \mathbb{W}(F) = \begin{cases} W^{\rm inc}(F) & \text{if } F \in \mathrm{SL}(2) \\ \infty & \text{if } F \notin \mathrm{SL}(2) \end{cases} \tag{2.10}$$

is polyconvex[3] another criterion which may be used for checking the polyconvexity of the energies defined on $\mathrm{SL}(n)$[4] is given in the Appendix of Neff et al. (2015b). A two-dimensional version may be formulated and proven exactly as in Neff et al. (2015b), as shown in the following:

**Proposition 2.5.** *Let $W : \mathrm{GL}^+(2) \to \mathbb{R}$ be an isotropic polyconvex function and let $g$ be a symmetric real-valued function defined on $\mathbb{R}^2_+$ such that, for all $F \in \mathrm{GL}^+(2)$, $W(F) = g(\lambda_1, \lambda_2)$, where $\lambda_1, \lambda_2$ are the singular values of $F$. Then there exists a convex function $\Psi : \mathbb{R}^3 \to \mathbb{R}_\infty := \mathbb{R} \cup \{+\infty\}$ with*

$$g(\lambda_1, \lambda_2) = \Psi(\lambda_1,\ \lambda_2, \lambda_1\lambda_2) \tag{2.11}$$

*for all $\lambda_1, \lambda_2 \in \mathbb{R}_+$.*

This necessary criteria was obtained in the Appendix of Neff et al. (2015b) and it is enough for the proof of the main result on the present paper, but we mention a related necessary and sufficient criteria recently given by Wiedemann and Peter in Wiedemann and Peter (2026), see also (Dacorogna, 2008, Thm. 5.43), in terms of the signed singular values defined as in the following for $n = 2$.

Every matrix $F \in \mathbb{R}^{2\times 2}$ admits a *signed singular value decomposition* of the form

$$F = R_1\,\mathrm{diag}(\nu)\,R_2\,, \tag{2.12}$$

where $R_1, R_2 \in \mathrm{SO}(2)$ and $\nu \in \mathbb{R}^2$. Here, $\mathrm{diag}(\nu) \in \mathbb{R}^{2\times 2}$ denotes the diagonal matrix with entries $\nu_1, \nu_2 \in \mathbb{R}$. Importantly, the entries in $\nu$ need not be positive, unlike singular values, hence the name *signed singular values*. The advantage over the classic singular value decomposition in the context of isotropic elasticity is the fact that we deal with SO(2) instead of O(2). The prize we pay is that the signed singular value decomposition is non-unique[5] in the sense that we may change the sign of both $\nu_1, \nu_2$ and absorb it into either $R_1$ or $R_2$, which then remain proper. In addition, we have the permutation invariance in $\nu_1, \nu_2$ analogous to (1.1). In light of these symmetries, we may state the following theorem:

**Theorem 2.6** (*Polyconvexity in terms of the signed singular values (Wiedemann and Peter, 2026)*). *An isotropic function $\mathbb{W} : \mathbb{R}^{2\times 2} \to \mathbb{R}_\infty$ is lower semicontinuous and polyconvex if and only if there exists a lower semicontinuous* **convex** *function $\Psi : \mathbb{R}^3 \to \mathbb{R}_\infty$, such that*

$$\mathbb{W}(F) = \Psi(\nu_1, \nu_2, \nu_1\nu_2) = \Psi(\nu_2, \nu_1, \nu_1\nu_2) = \Psi(-\nu_1, -\nu_2, \nu_1\nu_2), \tag{2.13}$$

*for all $F \in \mathbb{R}^{2\times 2}$ and for all choices of signed singular values $\nu \in \mathbb{R}^2$ of $F$.*

This result is closely related to the necessary and sufficient conditions for polyconvexity in 2D (Dacorogna and Koshigoe, 1993), there called *diagonal polyconvexity*.

### *2.2. Kirchhoff stress tensor and its principal components for functions on SL(n)*

#### *2.2.1. Kirchhoff stress tensor for functions on SL(n)*

We now recall some related notions of monotonicity. Let

$$\mathcal{M} := \{V \in \mathrm{Sym}^+(n) \mid \det V = 1\} \tag{2.14}$$

---

[2] Throughout this article, $\|X\|^2 = \langle X, X\rangle$ denotes the Frobenius tensor norm of $X \in \mathbb{R}^{n\times n}$, where $\langle X, Y\rangle = \mathrm{tr}(Y^T X)$ is the standard Euclidean scalar product on $\mathbb{R}^{n\times n}$. The identity tensor on $\mathbb{R}^{n\times n}$ will be denoted by $\mathbb{1}$, so that $\mathrm{tr}(X) = \langle X, \mathbb{1}\rangle$.

[3] In the sense of the definition given by (Ball (Ball, 1977)): A function $\mathbb{W} : \mathbb{R}^{2\times 2} \to \mathbb{R} \cup \{\infty\}$ is called polyconvex if there exists a convex function $P : \mathbb{R}^{2\times 2} \times \mathbb{R} \to \mathbb{R} \cup \{\infty\}$ such that $\mathbb{W}(F) = P(F, \det F)$ for all $F \in \mathbb{R}^{n\times n}$,

[4] A sufficient criteria is given by Ball:
**Theorem 2.4** ((Ball, 1977, page 367), 2D Sufficient Conditions for Polyconvexity of Isotropic Functions): Let $W(F) = \Psi(\lambda_1, \lambda_2, \lambda_1\lambda_2)$, where $\lambda_1, \lambda_2$ are the singular values of $F \in \mathrm{GL}^+(2)$, and
(a) $\Psi : \mathbb{R}^3_+ \to \mathbb{R}$ is convex,
(b) $\Psi(\widetilde{P}\,x, \delta) = \Psi(x, \delta)$ for all $\widetilde{P} \in \mathcal{P}_2$ (an element $\widetilde{P}$ of $\mathcal{P}_2$, acts on a vector $v \in \mathbb{R}^2$ by permuting its entries) and all $x \in \mathbb{R}^2_+$, $\delta \in \mathbb{R}_+$,
(c) $\Psi(x_1, x_2, \delta)$ is nondecreasing in each $x_i$, individually.
Then $W$ is polyconvex on $\mathrm{GL}^+(2)$.

[5] Some authors, e.g., Šilhavý (2000, Sect. 1) or (Dacorogna, 2008, Sect. 5.3.3), make an additional definition to render the signed singular decomposition unique, only to end up with the same ensuing symmetries for the energy function.

be the manifold of isochoric left stretch tensors. This manifold is naturally associated with $\mathrm{SL}(n)$: for $F \in \mathrm{SL}(n)$ with polar decomposition $F = R\,U = V\,R$ (see Neff et al. (2014)), one has $V \in \mathcal{M}$. Consider the set

$$\operatorname{dev}\,\mathrm{Sym}(n) := \{X \in \mathrm{Sym}(n) \mid \operatorname{tr} X = 0\}. \tag{2.15}$$

Let $W^{\rm inc} : \mathrm{SL}(n) \to \mathbb{R}$ be an isotropic and objective stored energy and $\mathcal{W}^{\rm inc} : \mathcal{M} \to \mathbb{R}$

$$W^{\rm inc}(F) = \mathcal{W}^{\rm inc}(V) \tag{2.16}$$

its representation in terms of the left stretch tensor $V$. The map

$$L : \mathcal{M} \to \operatorname{dev}\,\mathrm{Sym}(n), \qquad L(V) = \log V \tag{2.17}$$

is a global diffeomorphism with inverse $X \mapsto \exp X$. Thus, $\mathcal{M}$ can be identified with the linear space $\operatorname{dev}\mathrm{Sym}(n)$ via the logarithmic map. Define the log-reparametrized energy $\widehat{W}^{\rm inc} : \operatorname{dev}\,\mathrm{Sym}(n) \to \mathbb{R}$

$$\widehat{W}^{\rm inc}(X) := \mathcal{W}^{\rm inc}(\exp X), \quad \forall\, X \in \operatorname{dev}\,\mathrm{Sym}(n) \qquad \iff \qquad \widehat{W}^{\rm inc}(\log V) = \mathcal{W}^{\rm inc}(V) \quad \forall\, V \in \mathcal{M}. \tag{2.18}$$

Since $\operatorname{dev}\mathrm{Sym}(n)$ is a linear subspace of $\mathbb{R}^{n\times n}$, one has that the tangent space is given by

$$T_X\big(\operatorname{dev}\,\mathrm{Sym}(n)\big) = \operatorname{dev}\,\mathrm{Sym}(n) \qquad \forall\, X \in \operatorname{dev}\,\mathrm{Sym}(n). \tag{2.19}$$

We say that $\widehat{W}^{\rm inc}$ is differentiable at $X \in \operatorname{dev}\mathrm{Sym}(n)$ if there exists a linear mapping

$$\mathrm{D}_X \widehat{W}^{\rm inc}(X) : \operatorname{dev}\mathrm{Sym}(n) \to \mathbb{R} \tag{2.20}$$

such that

$$\mathrm{D}_X \widehat{W}^{\rm inc}(X).H = \lim_{t\to 0} \frac{\widehat{W}^{\rm inc}(X+tH) - \widehat{W}^{\rm inc}(X)}{t} \qquad \forall\, H \in \operatorname{dev}\mathrm{Sym}(n). \tag{2.21}$$

Identifying $(\operatorname{dev}\mathrm{Sym}(n))^*$ with $\operatorname{dev}\mathrm{Sym}(n)$ via the Frobenius inner product, by the Riesz representation theorem there exists a unique $\widehat{\boldsymbol{\tau}}(X) \in \operatorname{dev}\mathrm{Sym}(n)$ such that

$$\langle \widehat{\boldsymbol{\tau}}(X), H\rangle = \mathrm{D}_X \widehat{W}^{\rm inc}(X).H \qquad \forall\, H \in \operatorname{dev}\mathrm{Sym}(n). \tag{2.22}$$

Equivalently, for $V \in \mathcal{M}$ and $X = \log V$, one has

$$\langle \widehat{\boldsymbol{\tau}}(\log V), H\rangle = \mathrm{D}_{\log V} \widehat{W}^{\rm inc}(\log V).H \qquad \forall\, H \in \operatorname{dev}\mathrm{Sym}(n). \tag{2.23}$$

We refer to $\widehat{\boldsymbol{\tau}} : \operatorname{dev}\mathrm{Sym}(n) \to \operatorname{dev}\mathrm{Sym}(n)$ as the Kirchhoff stress tensor for energies defined on $\mathrm{SL}(n)$. This is in accord with the Richter-formula (see, e.g., Richter (1949, 1948))

$$\widehat{\boldsymbol{\tau}}(\log V) = \mathrm{D}_{\log V} \widehat{W}^{\rm inc}(\log V). \tag{2.24}$$

Let us notice that we are working with functions defined on $\operatorname{dev}\mathrm{Sym}(n)$. If we assume that the energy $\widehat{W}^{\rm inc} : \operatorname{dev}\mathrm{Sym}(n) \to \mathbb{R}$ is the restriction of a function $\widehat{W} : \mathrm{Sym}(n) \to \mathbb{R}$, then the derivative can be computed in the ambient space $\mathrm{Sym}(n)$ by introducing a Lagrange multiplier associated with the constraint $\operatorname{tr} X = 0$ and we have[6]

$$\widehat{\boldsymbol{\tau}}(X) = \operatorname{dev}\big(\mathrm{D}_X \widehat{W}(X)\big) \qquad \iff \qquad \widehat{\boldsymbol{\tau}}(\log V) = \operatorname{dev}\big(\mathrm{D}_{\log V} \widehat{W}(\log V)\big) \quad \forall\, V \in \mathrm{SL}(n). \tag{2.25}$$

[6] More precisely, consider the augmented functional $\mathcal{L}(X,p) := \widehat{W}(X) + p\operatorname{tr} X$, where $p \in \mathbb{R}$ is a scalar Lagrange multiplier. Then, for every $H \in \mathrm{Sym}(n)$, one has $\mathrm{D}_X\mathcal{L}(X,p).h = \langle \mathrm{D}_X\widehat{W}(X), H\rangle + p\operatorname{tr} H$. Restricting to admissible variations $H \in \operatorname{dev}\mathrm{Sym}(n)$ (i.e. $\operatorname{tr} H = 0$), we recover the derivative of $\widehat{W}^{\rm inc}$: $\mathrm{D}_X\widehat{W}^{\rm inc}(X).h = \langle \mathrm{D}_X\widehat{W}(X), H\rangle = \langle (\mathrm{D}_X\widehat{W}(X) + p\,\mathbb{1}), H\rangle$. Therefore, the Kirchhoff stress tensor admits the representation $\widehat{\boldsymbol{\tau}}(X) = \mathrm{D}_X\widehat{W}(X) + p(X)\,\mathbb{1}$, $X \in \operatorname{dev}\mathrm{Sym}(n)$, where $p(X)$ is the Lagrange multiplier associated with the constraint $\operatorname{tr} X = 0$. Imposing the condition $\widehat{\boldsymbol{\tau}}(X) \in \operatorname{dev}\mathrm{Sym}(n)$ yields $p(X) = -\frac{1}{n}\operatorname{tr}\big(\mathrm{D}_X\widehat{W}(X)\big)$, and thus $\widehat{\boldsymbol{\tau}}(X) = \operatorname{dev}\big(\mathrm{D}_X\widehat{W}(X)\big)$. The scalar field $p(X)$ represents the indeterminate spherical part of the stress tensor, and is completely analogous to the pressure in incompressible elasticity.

### 2.2.2. The vector of the principal components of the Kirchhoff stress tensor on SL($n$)

Similarly, for $\ell_i = \log\lambda_i$, $i = 1,2,\ldots,n$ consider

$$\begin{aligned} \Pi &:= \left\{ \ell \in \mathbb{R}^n \;\middle|\; \sum_{i=1}^n \ell_i = 0 \right\} \\ &= \left\{ (\log\lambda_1, \log\lambda_2, \ldots, \log\lambda_n) \in \mathbb{R}^n \;\middle|\; \sum_{i=1}^n \log\lambda_i = 0 \right\}. \end{aligned} \tag{2.26}$$

Then $\Pi$ is a linear subspace of $\mathbb{R}^n$, and $T_\ell \Pi = \Pi$, $\forall \ell \in \Pi$. Let $\widehat{g}^{\rm inc} : \Pi \to \mathbb{R}$ be the differentiable function such that the isotropic energy satisfies $\widehat{W}^{\rm inc}(X) = \widehat{g}(\ell)$, $X \in \operatorname{dev}\mathrm{Sym}(n)$, $\ell \in \Pi$. The intrinsic derivative at $\ell$ is the linear map $\mathrm{D}_\ell \widehat{g}^{\rm inc}(\ell) : \Pi \to \mathbb{R}$ defined by

$$\mathrm{D}_\ell \widehat{g}^{\rm inc}(\ell).h = \lim_{t\to 0} \frac{\widehat{g}^{\rm inc}(\ell + t\,h) - \widehat{g}^{\rm inc}(\ell)}{t} \qquad \forall\, h \in \Pi. \tag{2.27}$$

Identifying $\Pi^*$ with $\Pi$ via the Euclidean inner product, by the Riesz representation theorem there exists a unique vector $\widehat{\tau}(\ell) \in \Pi$ such that

$$\langle \widehat{\tau}(\ell), h\rangle = \mathrm{D}_\ell \widehat{g}^{\rm inc}(\ell).h \qquad \forall\, h \in \Pi. \tag{2.28}$$

For any vector $w \in \mathbb{R}^n$, we define its "deviatoric part" by $\operatorname{dev}(w) := w - \frac{1}{n}\Big(\sum_{i=1}^n w_i\Big)\mathbf{1}$, with $\mathbf{1} = (1,\ldots,1)$. This is the orthogonal projection of $w$ onto $\Pi$. Let $\widehat{g}^{\rm inc} : \Pi \to \mathbb{R}$ be a differentiable function, and assume that it is the restriction of a smooth function $g : \mathbb{R}^n \to \mathbb{R}$. To compute derivatives intrinsically on $\Pi$, we introduce a Lagrange multiplier associated with the constraint $\sum_i \ell_i = 0$ and we have[7]

$$\widehat{\tau}(\ell) = \operatorname{dev}\big(\mathrm{D}_\ell g(\ell)\big). \tag{2.29}$$

Thus, the intrinsic gradient of $\widehat{g}^{\rm inc}$ on $\Pi$ is the orthogonal projection of the ambient gradient $\mathrm{D}_\ell g(\ell)$ onto $\Pi$. We call $\widehat{\tau}$ the vector of the principal components of the Kirchhoff stress tensor. However, in order to justify this definition, we have to show that indeed the thus defined vector $\widehat{\tau}$ is the vector containing the eigenvalues of the Kirchhoff stress tensor $\widehat{\boldsymbol{\tau}}$.

### 2.2.3. From the Kirchhoff stress tensor to its principal components and vice-versa

We show now the relation between the intrinsic gradient of $\widehat{g}^{\rm inc}$ and the intrinsic derivative of $\widehat{W}^{\rm inc}$. Although the main idea of the proof remains largely unchanged, the incompressible case requires a separate argument since the admissible tensors no longer vary in all directions but only along a constrained subset corresponding to volume preservation. Consequently, monotonicity must be verified only along directions compatible with this constraint, and results from the compressible setting do not automatically transfer. Moreover, one cannot rely on properties of possible extensions outside the incompressible setting, since different extensions may exhibit different monotonicity properties, with some being monotone and others not. To the best of our knowledge, such a proof has not been carried out explicitly in the literature, and therefore cannot simply be cited.

[7] For this, let us define the augmented functional $\mathcal{L}(\ell,p) = g(\ell) + p\sum_{i=1}^n \ell_i$, $(\ell,p) \in \mathbb{R}^n \times \mathbb{R}$. Then, for every $h \in \mathbb{R}^n$, we have $\mathrm{D}_\ell\mathcal{L}(\ell,p).h = \langle \mathrm{D}_\ell g(\ell), h\rangle + p\sum_{i=1}^n h_i$. Restricting to admissible variations $h \in \Pi$ (i.e. $\sum_i h_i = 0$), we obtain $\mathrm{D}_\ell \widehat{g}^{\rm inc}(\ell).h = \langle \mathrm{D}_\ell g(\ell), h\rangle = \langle \mathrm{D}_\ell g(\ell) + p(\ell)\mathbf{1}, h\rangle$. By the Riesz representation theorem, there exists a unique vector $\widehat{\tau}(\ell) \in \Pi$ such that $\langle \widehat{\tau}(\ell), h\rangle = \mathrm{D}_\ell \widehat{g}^{\rm inc}(\ell).h \quad \forall\, h \in \Pi$. Hence if we take $\widehat{\tau}(\ell) = \mathrm{D}_\ell g(\ell) + p(\ell)\mathbf{1}$ then, imposing the constraint $\widehat{\tau}(\ell) \in \Pi$ yields $0 = \sum_{i=1}^n \widehat{\tau}_i(\ell) = \sum_{i=1}^n \partial_i g(\ell) + n\,p(\ell)$, and therefore $p(\ell) = -\frac{1}{n}\sum_{i=1}^n \partial_i g(\ell)$. Substituting back, we obtain the intrinsic gradient $\widehat{\tau}(\ell) = \mathrm{D}_\ell g(\ell) - \frac{1}{n}\Big(\sum_{i=1}^n \partial_i g(\ell)\Big)\mathbf{1} = \operatorname{dev}\big(\mathrm{D}_\ell g(\ell)\big)$.

**Theorem 2.7.** *Let $\widehat{W}^{\rm inc} : \operatorname{dev}\operatorname{Sym}(n) \to \mathbb{R}$ be an isotropic differentiable function, and let*

$$X = Q\operatorname{diag}(\ell)\, Q^T, \qquad \ell = (\ell_1, \dots, \ell_n) \in \Pi, \qquad Q \in \mathrm{O}(n), \tag{2.30}$$

*be the spectral decomposition of $X$. Assume that $\widehat{W}^{\rm inc}$ admits the spectral representation*

$$\widehat{W}^{\rm inc}(X) = \widehat{g}^{\rm inc}(\ell) \qquad \Longleftrightarrow \qquad \widehat{W}^{\rm inc}(\log V) = \widehat{g}^{\rm inc}(\log\lambda_1, \log\lambda_2, \dots, \log\lambda_n) \tag{2.31}$$

*for some symmetric differentiable function $\widehat{g}^{\rm inc} : \Pi \to \mathbb{R}$. Then*

$$\widehat{\boldsymbol{\tau}}(X) = Q\operatorname{diag}(\widehat{\tau}(\ell))\, Q^T. \tag{2.32}$$

*In particular, the components of $\widehat{\tau}(\ell)$ are exactly the principal values of the Kirchhoff stress tensor $\widehat{\boldsymbol{\tau}}(X)$.*

**Proof.** Let $h \in \Pi$ and define $H := Q\operatorname{diag}(h)\, Q^T$. Since $\sum_{i=1}^n h_i = 0$, we have $\operatorname{tr} H = \operatorname{tr}(\operatorname{diag}(h)) = \sum_{i=1}^n h_i = 0$, hence $H \in \operatorname{dev}\operatorname{Sym}(n)$. Now consider the curve $X(t) := X + tH = Q\operatorname{diag}(\ell + t\, h)\, Q^T$. Because $\widehat{W}^{\rm inc}$ is isotropic and admits the spectral representation $\widehat{W}^{\rm inc}(Q\operatorname{diag}(w)\, Q^T) = \widehat{g}^{\rm inc}(w)$, we obtain $\widehat{W}^{\rm inc}(X + tH) = \widehat{g}^{\rm inc}(\ell + t\, h)$. Differentiating at $t = 0$ yields $\mathrm{D}\widehat{W}^{\rm inc}(X).H = \mathrm{D}\widehat{g}^{\rm inc}(\ell).h$

Using the definitions of $\widehat{\boldsymbol{\tau}}(X)$ and $\widehat{\tau}(\ell)$ through the Riesz representation theorem, we therefore get

$$\langle \widehat{\boldsymbol{\tau}}(X), H\rangle = \mathrm{D}\widehat{W}^{\rm inc}(X).H = \mathrm{D}\widehat{g}^{\rm inc}(\ell).h = \langle \widehat{\tau}(\ell), h\rangle. \tag{2.33}$$

Since $H = Q\operatorname{diag}(h)\, Q^T$, the left-hand side can be rewritten as

$$\langle \widehat{\boldsymbol{\tau}}(X), H\rangle = \langle \widehat{\boldsymbol{\tau}}(X), Q\operatorname{diag}(h)\, Q^T\rangle = \langle Q^T\widehat{\boldsymbol{\tau}}(X)\, Q, \operatorname{diag}(h)\rangle \tag{2.34}$$

Since $\langle Q^T\widehat{\boldsymbol{\tau}}(X)\, Q, \operatorname{diag}(h)\rangle = \langle \widehat{\tau}(\ell), h\rangle$ and
$\langle Q^T\widehat{\boldsymbol{\tau}}(X)\, Q, \operatorname{diag}(h)\rangle = \sum_{i=1}^n (Q^T\widehat{\boldsymbol{\tau}}(X)\, Q)_{ii} h_i$, it follows that

$$\sum_{i=1}^n \left[(Q^T\widehat{\boldsymbol{\tau}}(X)\, Q)_{ii} - \widehat{\tau}_i(\ell)\right]\, h_i = 0 \qquad \forall\, h \in \Pi. \tag{2.35}$$

Therefore the vector

$$\big((Q^T\widehat{\boldsymbol{\tau}}(X)\, Q)_{11} - \widehat{\tau}_1(\ell), \dots, (Q^T\widehat{\boldsymbol{\tau}}(X)\, Q)_{nn} - \widehat{\tau}_n(\ell)\big) \tag{2.36}$$

is orthogonal to $\Pi$. Since $\Pi^\perp = \operatorname{span}\{(1, \dots, 1)\}$, there exists $c \in \mathbb{R}$ such that

$$(Q^T\widehat{\boldsymbol{\tau}}(X)\, Q)_{ii} - \widehat{\tau}_i(\ell) = c \qquad \forall i = 1, \dots, n. \tag{2.37}$$

Now both $Q^T\widehat{\boldsymbol{\tau}}(X)\, Q$ and $\operatorname{diag}(\widehat{\tau}(\ell))$ are trace-free, because $(Q^T\widehat{\boldsymbol{\tau}}(X)\, Q) \in \operatorname{dev}\operatorname{Sym}(n)$, $\widehat{\tau}(\ell) \in \Pi$. Thus

$$0 = \operatorname{tr}(Q^T\widehat{\boldsymbol{\tau}}(X)\, Q) - \operatorname{tr}(\operatorname{diag}(\widehat{\tau}(\ell))) = \sum_{i=1}^n \big((Q^T\widehat{\boldsymbol{\tau}}(X)\, Q)_{ii} - \widehat{\tau}_i(\ell)\big) = n\, c, \tag{2.38}$$

and therefore $c = 0$. Consequently,

$$(Q^T\widehat{\boldsymbol{\tau}}(X)\, Q)_{ii} = \widehat{\tau}_i(\ell) \qquad \forall i = 1, \dots, n. \tag{2.39}$$

It remains to show that the off-diagonal entries of $(Q^T\widehat{\boldsymbol{\tau}}(X)\, Q)$ vanish. Since $\widehat{W}^{\rm inc}$ is isotropic, one has $\widehat{W}^{\rm inc}(QXQ^T) = \widehat{W}^{\rm inc}(X)\forall\, Q \in \mathrm{O}(n)$. Differentiating this identity with respect to $X$ yields

$$\widehat{\boldsymbol{\tau}}(QXQ^T) = Q\,\widehat{\boldsymbol{\tau}}(X)\, Q^T. \tag{2.40}$$

Let $Q$ be such that $X = Q\operatorname{diag}(\ell_1, \dots, \ell_n)\, Q^T$. Then

$$Q^T\widehat{\boldsymbol{\tau}}(X)\, Q = \widehat{\boldsymbol{\tau}}(\operatorname{diag}(\ell_1, \dots, \ell_n)). \tag{2.41}$$

Since $\operatorname{diag}(\ell_1, \dots, \ell_n)$ is diagonal, it follows that $S := Q^T\widehat{\boldsymbol{\tau}}(X)\, Q$ is a function of a diagonal matrix. Since $\widehat{W}^{\rm inc}$ is isotropic, one has $\widehat{\boldsymbol{\tau}}(PDP^T) = P\widehat{\boldsymbol{\tau}}(D)P^T, \forall P \in O(n)$. In particular, if $PDP^T = D$, then $S = PSP^T$. Choosing $P = \operatorname{diag}(\varepsilon_1, \dots, \varepsilon_n)$ with $\varepsilon_i = \pm 1$, we obtain $S_{ij} = \varepsilon_i\varepsilon_j S_{ij}$. For $i \neq j$, selecting $\varepsilon_i = -1$, $\varepsilon_j = 1$ yields $S_{ij} = 0$. Thus $S := Q^T\widehat{\boldsymbol{\tau}}(X)\, Q$ is diagonal having on the diagonal the vector $\widehat{\tau}(\ell)$. ■

Therefore $\widehat{\boldsymbol{\tau}} : \operatorname{dev}\operatorname{Sym}(n) \to \operatorname{Sym}(n)$ is an isotropic tensor function satisfying

$$\widehat{\boldsymbol{\tau}}(Q^T\operatorname{diag}(\ell_1, \ell_2, \dots, \ell_n)\, Q) = Q^T\widehat{\boldsymbol{\tau}}(\operatorname{diag}(\ell_1, \ell_2, \dots, \ell_n))\, Q \tag{2.42}$$

and

$$\widehat{\boldsymbol{\tau}}(\underbrace{Q^T\operatorname{diag}(\ell_1, \ell_2, \dots, \ell_n)\, Q}_{X \in \operatorname{dev}\operatorname{Sym}(n)}) = \underbrace{Q^T\operatorname{diag}(\widehat{\tau}(\ell_1, \ell_2, \dots, \ell_n))\, Q}_{\widehat{\boldsymbol{\tau}}(X) \in \operatorname{Sym}(n)} \quad \forall\, Q \in \mathrm{O}(n) \tag{2.43}$$

with a vector-function $\widehat{\tau} = (\widehat{\tau}_1, \widehat{\tau}_2, \dots, \widehat{\tau}_n) : \Pi \to \mathbb{R}^n$,which fulfills

$$\widehat{\tau}_i(\ell_{\pi(1)}, \ell_{\pi(2)}, \dots, \ell_{\pi(n)}) = \widehat{\tau}_{\pi(i)}(\ell_1, \ell_2, \dots, \ell_n) \tag{2.44}$$

for any permutation $\pi : \{1, 2, \dots, n\} \to \{1, 2, \dots, n\}$. The permutation symmetry implies and it is implied by isotropy.

### 2.3. Hilbert-monotonicity and Hill's inequality

The results obtained in the previous subsection clarify why the following definitions are well formulated.

**Definition 2.8.** The Kirchhoff stress tensor $\widehat{\boldsymbol{\tau}} : \operatorname{dev}\operatorname{Sym}(n) \to \operatorname{dev}\operatorname{Sym}(n)$ is called **strictly Hilbert-monotone** if

$$\begin{aligned}
&\langle \widehat{\boldsymbol{\tau}}(X) - \widehat{\boldsymbol{\tau}}(\overline{X}),\, X - \overline{X}\rangle_{\mathbb{R}^{n\times n}} > 0 \quad \forall\, X \neq \overline{X} \in \operatorname{dev}\operatorname{Sym}(n)\\
\Longleftrightarrow \quad &\langle \widehat{\boldsymbol{\tau}}(\log V) - \widehat{\boldsymbol{\tau}}(\log \overline{V}), \log V - \log \overline{V}\rangle_{\mathbb{R}^{n\times n}} > 0 \quad \forall\, V \neq \overline{V} \in \operatorname{Sym}^+(n).
\end{aligned} \tag{2.45}$$

We refer to this inequality as **strict Hilbert-space matrix-monotonicity** of the tensor function $\widehat{\boldsymbol{\tau}}$. The Kirchhoff stress tensor $\widehat{\boldsymbol{\tau}} : \operatorname{dev}\operatorname{Sym}(n) \to \operatorname{dev}\operatorname{Sym}(n)$ is called **Hilbert-monotone** if

$$\langle \widehat{\boldsymbol{\tau}}(X) - \widehat{\boldsymbol{\tau}}(\overline{X}),\, X - \overline{X}\rangle_{\mathbb{R}^{n\times n}} \geq 0 \quad \forall\, X, \overline{X} \in \operatorname{dev}\operatorname{Sym}(n). \tag{2.46}$$

**Definition 2.9.** The vector function $\widehat{\tau} : \Pi \to \Pi$ is **strictly vector monotone** if

$$\langle \widehat{\tau}(\ell) - \widehat{\tau}(\overline{\ell}),\, \ell - \overline{\ell}\rangle_{\mathbb{R}^n} > 0 \quad \forall\, \ell \neq \overline{\ell} \in \Pi, \tag{2.47}$$

and it is **vector monotone** if

$$\langle \widehat{\tau}(\ell) - \widehat{\tau}(\overline{\ell}),\, \ell - \overline{\ell}\rangle_{\mathbb{R}^n} \geq 0 \quad \forall\, \ell, \overline{\ell} \in \Pi, \tag{2.48}$$

**Definition 2.10.** If $\widehat{\boldsymbol{\tau}} : \operatorname{dev}\operatorname{Sym}(n) \to \operatorname{dev}\operatorname{Sym}(n)$ is a continuously differentiable tensor function, then it is called **strongly Hilbert-monotone** if

$$\langle \mathrm{D}_X\widehat{\boldsymbol{\tau}}(X).H, H\rangle > 0 \qquad \text{for all} \quad X \in \operatorname{dev}\operatorname{Sym}(n),\ H \in \operatorname{dev}\operatorname{Sym}(n). \tag{2.49}$$

**Definition 2.11.** If $\widehat{\tau} : \Pi \to \Pi$ is a continuously differentiable vector function, then it is called **strongly vector monotone** if

$$\langle \mathrm{D}_\ell\widehat{\tau}(\ell).h, h\rangle > 0 \qquad \text{for all} \quad \ell \in \Pi,\ h \in \Pi. \tag{2.50}$$

In a forthcoming paper (Martin et al., 2026), we discuss a result concerning matrix monotonicity versus vector monotonicity in the context of compressible isotropic elastic materials, see also Ogden's work (Ogden, 1983, last page in the Appendix) and Hill's seminal contributions (Hill, 1968b,a, 1970).

In Appendix A, we give some milestones of the proof of the following result to assure the reader that the conclusions regarding the matrix monotonicity vs. vector monotonicity of compressible isotropic elastic materials hold true for incompressible materials, too.

**Theorem 2.12.** *The Kirchhoff tensor $\widehat{\boldsymbol{\tau}} : \operatorname{dev}\operatorname{Sym}(n) \to \operatorname{dev}\operatorname{Sym}(n)$ is (strictly/strongly) matrix monotone if and only if $\widehat{\tau} : \Pi \to \Pi$ is (strictly/strongly) vector-monotone.*

**Proposition 2.13** (*Equivalence of vector monotonicity and differential monotonicity*). *Let* $\hat{\tau} : \Pi \to \Pi$, *where* $\Pi := \{\ell \in \mathbb{R}^n \mid \sum_{i=1}^n \ell_i = 0\}$, *be a continuously differentiable mapping. Then the following are equivalent:*

*(i)* $\hat{\tau}$ *is vector monotone (strict monotone).*
*(ii) The Jacobian of* $\hat{\tau}$ *is positive semidefinite (positive definite), i.e.*

$$\langle \mathrm{D}_\ell \hat{\tau}(\ell).\, h,\, h\rangle \geq 0 \qquad \forall\, \ell \in \Pi,\ h \in \Pi. \tag{2.51}$$

**Proof.** (i) ⇒ (ii): Let $\ell \in \Pi$ and $h \in \Pi$. For $t > 0$, apply monotonicity with $\bar{\ell} = \ell$ and $\ell_t = \ell + th$: $\langle \hat{\tau}(\ell + t\,h) - \hat{\tau}(\ell), t\,h\rangle \geq 0$. Dividing by $t^2$ yields $\left\langle \frac{\hat{\tau}(\ell+t\,h)-\hat{\tau}(\ell)}{t}, h\right\rangle \geq 0$. Passing to the limit as $t \to 0$ gives $\langle \mathrm{D}_\ell \hat{\tau}(\ell).\, h,\, h\rangle \geq 0$.

(ii) ⇒ (i): Let $\ell, \bar{\ell} \in \Pi$ and define the segment $\gamma(t) := \bar{\ell} + t(\ell - \bar{\ell})$, $t \in [0,1]$. Set $\phi(t) := \langle \hat{\tau}(\gamma(t)), \ell - \bar{\ell}\rangle$. Then $\phi'(t) = \langle \mathrm{D}\hat{\tau}(\gamma(t))[\ell - \bar{\ell}], \ell - \bar{\ell}\rangle \geq 0$. Hence $\phi$ is non-decreasing, and therefore $\phi(1) - \phi(0) \geq 0$. This gives $\langle \hat{\tau}(\ell) - \hat{\tau}(\bar{\ell}), \ell - \bar{\ell}\rangle \geq 0$. The strict case follows analogously. ■

## 3. Polyconvexity implies the weak Hill inequality on SL(2) for objective and isotropic energies

In the following two subsections we give alternative proofs of Theorem 1.2.

### 3.1. A proof based on Theorem 2.3

Since $W^{\mathrm{inc}}$ is objective and isotropic, there exists a function $\phi : [0,\infty) \to \mathbb{R}$ such that

$$W^{\mathrm{inc}}(F) = \phi\left(\sqrt{\|F\|^2 - 2}\right) = \phi\left(\lambda_{\max}(F) - \frac{1}{\lambda_{\max}(F)}\right) \quad \forall\, F \in \mathrm{SL}(2). \tag{3.1}$$

Because $W^{\mathrm{inc}}$ is polyconvex, Theorem 2.3 implies that $\phi$ is convex and nondecreasing on $[0,\infty)$. Now let $X \in \operatorname{dev}\operatorname{Sym}(2)$. Then

$$X = Q\,\mathrm{diag}(t, -t)\,Q^T \qquad \text{for some } Q \in O(2),\ t \in \mathbb{R}. \tag{3.2}$$

Hence $\exp(X) = Q\,\mathrm{diag}(e^t, e^{-t})\,Q^T$, and therefore

$$\widehat{W}^{\mathrm{inc}}(X) = \phi\left(e^{|t|} - e^{-|t|}\right) = \phi(2\sinh|t|)\,. \tag{3.3}$$

Define

$$\psi : \mathbb{R} \to \mathbb{R}, \qquad \psi(t) := \phi(2\sinh|t|). \tag{3.4}$$

Then

$$\widehat{W}^{\mathrm{inc}}(Q\,\mathrm{diag}(t,-t)\,Q^T) = \psi(t). \tag{3.5}$$

The function $\alpha(t) := 2\sinh|t|$ is convex on $\mathbb{R}$ and takes values in $[0,\infty)$. Since $\phi$ is convex and nondecreasing on $[0,\infty)$, the composition $\psi = \phi\circ\alpha$ is convex on $\mathbb{R}$. Now observe that the reduced spectral function

$$\widehat{W}^{\mathrm{inc}}(X) = \hat{g}^{\mathrm{inc}}(\ell_1, \ell_2) \tag{3.6}$$

on $\Pi := \{(\ell_1,\ell_2) \in \mathbb{R}^2 \mid \ell_1 + \ell_2 = 0\}$ is given by

$$\hat{g}^{\mathrm{inc}}(t, -t) := \psi(t). \tag{3.7}$$

Since $\psi$ is convex on $\mathbb{R}$, the map $\hat{g}^{\mathrm{inc}}$ is convex on the one-dimensional affine space $\Pi$. Because $\widehat{W}^{\mathrm{inc}}$ is differentiable, its intrinsic gradient on $\Pi$, denoted by $\hat{\tau} : \Pi \to \Pi$, exists and is monotone:

$$\langle \hat{\tau}(\ell) - \hat{\tau}(h), \ell - h\rangle \geq 0 \quad \forall\, \ell, h \in \Pi. \tag{3.8}$$

Indeed, this is the standard monotonicity of the gradient of a convex differentiable function on a Euclidean space. By isotropy, the Kirchhoff stress in logarithmic variables has the spectral representation

$$\hat{\boldsymbol{\tau}}(X) = Q\,\mathrm{diag}(\hat{\tau}(\ell))\,Q^T \qquad \text{whenever } X = Q\,\mathrm{diag}(\ell)\,Q^T,\quad \ell \in \Pi. \tag{3.9}$$

Therefore, by the spectral monotonicity result from the previous section, see Theorem 2.12,

$$\langle \hat{\boldsymbol{\tau}}(X_1) - \hat{\boldsymbol{\tau}}(X_2),\, X_1 - X_2\rangle \geq 0 \quad \forall\, X_1, X_2 \in \operatorname{dev}\operatorname{Sym}(2). \tag{3.10}$$

Finally, taking $X_i = \log V_i$ with $X_i \in \mathrm{Sym}^+(2) \cap \mathrm{SL}(2)$ and using

$$\hat{\boldsymbol{\tau}}(\log V_i) = \boldsymbol{\tau}(\ell_i), \tag{3.11}$$

we obtain

$$\langle \boldsymbol{\tau}(V_1) - \boldsymbol{\tau}(V_2),\, \log V_1 - \log V_2\rangle \geq 0. \tag{3.12}$$

This is the (**weak**) Hill inequality and the proof is complete.

### 3.2. An alternative proof based on Theorem 2.6

Consider the extension $\mathbb{W}$ of $W^{\mathrm{inc}}(F)$ given by (2.10). Then the polyconvexity of $W^{\mathrm{inc}}(F)$ on SL(2) implies the existence of a lower semicontinuous convex function $\Psi : \mathbb{R}^{n\times n} \to \mathbb{R}_\infty$ by virtue of Theorem 2.6 such that

$$\mathbb{W}(F) = \Psi(\nu_1, \nu_2, \nu_1\nu_2). \tag{3.13}$$

Let $\vartheta = \vartheta(\nu_1, \nu_2)$ be defined as the slice

$$\vartheta(\nu_1, \nu_2) := \Psi(\nu_1, \nu_2, 1). \tag{3.14}$$

Then $\vartheta$ is convex on $\mathbb{R}^2$ due to the convexity on $\Psi$ on $\mathbb{R}^3$. Moreover, by isotropy,

$$\vartheta(\nu_1, \nu_2) = \vartheta(\nu_2, \nu_1), \tag{3.15}$$

and by invariance under pairwise reflection, see (2.13), we have

$$\vartheta(\nu_1, \nu_2) = \vartheta(-\nu_1, -\nu_2). \tag{3.16}$$

Now define

$$\widetilde{\vartheta}(x, y) := \vartheta\left(\frac{x+y}{2}, \frac{x-y}{2}\right) \quad \text{such that} \quad \vartheta(\nu_1, \nu_2) = \widetilde{\vartheta}(\nu_1 + \nu_2, \nu_1 - \nu_2). \tag{3.17}$$

Since the map

$$(x, y) \mapsto \left(\frac{x+y}{2}, \frac{x-y}{2}\right) \tag{3.18}$$

is linear and $\vartheta$ is convex, it follows that $\widetilde{\vartheta}$ is convex. The permutation-invariance of $\vartheta$ from (3.15) implies

$$\widetilde{\vartheta}(x, y) = \widetilde{\vartheta}(x, -y), \tag{3.19}$$

while the invariance under pairwise reflection from (3.16) implies

$$\widetilde{\vartheta}(x, y) = \widetilde{\vartheta}(-x, -y) = \widetilde{\vartheta}(-x, y). \tag{3.20}$$

Hence, $\widetilde{\vartheta}$ is even in each argument. In combination with convexity, it follows that

$$\begin{aligned} &\widetilde{\vartheta}(x_2, y) > \widetilde{\vartheta}(x_1, y) \quad \forall\, x_2 > x_1 \geq 0 \qquad \text{and} \\ &\widetilde{\vartheta}(x, y_2) > \widetilde{\vartheta}(x, y_1) \quad \forall\, y_2 > y_1 \geq 0, \end{aligned} \tag{3.21}$$

i.e., $\widetilde{\vartheta}$ is also strictly increasing[8] in each argument. Now observe that the reduced spectral function

$$\widehat{W}^{\mathrm{inc}}(X) = \hat{g}^{\mathrm{inc}}(\ell_1, \ell_2) \tag{3.22}$$

on $\Pi = \{(\ell_1,\ell_2) \in \mathbb{R}^2 \mid \ell_1 + \ell_2 = 0\}$ is given by

$$\hat{g}^{\mathrm{inc}}(\ell, -\ell) = \vartheta(e^\ell, e^{-\ell}) = \widetilde{\vartheta}(\cosh\ell, \sinh\ell) = \widetilde{\vartheta}(\cosh\ell, |\sinh\ell|). \tag{3.23}$$

Since the functions $\ell \mapsto \cosh\ell$ and $\ell \mapsto |\sinh\ell|$ are strictly convex and $\widetilde{\vartheta}$ is convex and strictly increasing in each variable on $[0,\infty)$, it follows that the map $\hat{g}^{\mathrm{inc}}$ is strictly convex on the one-dimensional affine space $\Pi$. The rest of the proof is similar as in the proof given in the previous section.

[8] The strictness of the implied monotonicity hinges on the assumption that $\widetilde{\vartheta}$ is non-constant in a region around the origin. Otherwise, $\widetilde{\vartheta}$ is only increasing.

## 4. Conclusion and outlook

In this contribution, we provided two proofs showing that polyconvexity (rank-one convexity) implies Hill's inequality and, in turn, true-stress–true-strain monotonicity (TSTS–M$^+$) in the two-dimensional incompressible case. In the process, we also revisited the transition from tensor arguments to principal-value representations for questions of convexity and monotonicity of isotropic functions.

Concerning the incompressible three-dimensional case, we know that Legendre–Hadamard ellipticity does not imply Hill's inequality. For instance, $W(F) = \sqrt{\mathrm{tr}(F^T F) - 3 + \alpha}$ for all $\alpha \in (0, 3 - 15 \times 2^{-7/3})$ is rank-one convex by virtue of Dunn et al. (2003), but does not satisfy Hill's inequality since in uniaxial tension the Cauchy stress response is not monotone in compression. Recently, Klein et al. (2026) also identified an incompressible polyconvex energy function that does not satisfy Hill's inequality.

The implication established in the present paper is essentially tied to the special structure of the two-dimensional incompressible setting. In SL(2), incompressibility and isotropy reduce the dependence on the principal stretches to a single independent scalar variable, and the available characterizations of rank-one convexity and polyconvexity can therefore be directly related to convexity in the logarithmic stretch variable. No analogous general reduction is available in three dimensions, where the incompressibility constraint leaves two independent principal stretches and the structure of polyconvexity is considerably more involved. Moreover, the implication "polyconvexity implies Hill's inequality" does not hold for general incompressible energies in three dimensions.

An intermediate open question is to find an energy defined on SL(3) that is Legendre–Hadamard elliptic but not polyconvex. The difficulty here is that necessary and sufficient conditions for polyconvexity only postulate the existence of a particular higher-dimensional function, but provide little to show lack thereof. A possible avenue would be to prove lack of quasi-convexity which is made difficult by the closed-form construction of an isochoric, heterogeneous perturbation field over which one could volume-average. We nonetheless guess that the construction of such an example is possible in the three-dimensional case, while in two dimensions all conditions coincide, as seen here.

Another open question is the relation between incompressible polyconvexity and Hill's inequality in 3D for a subclass of energy functions of Valanis–Landel type (Valanis and Landel, 1967), where the reduced functional space might indeed entail an implication similar to the general two-dimensional case discussed here.

### CRediT authorship contribution statement

**Ionel-Dumitrel Ghiba:** Writing – review & editing, Writing – original draft, Formal analysis, Conceptualization. **Maximilian P. Wollner:** Writing – review & editing, Writing – original draft, Formal analysis, Conceptualization. **Patrizio Neff:** Writing – review & editing, Writing – original draft, Formal analysis, Conceptualization.

### Declaration of competing interest

The authors declare that they have no known competing financial interests or personal relationships that could have appeared to influence the work reported in this paper.

### Acknowledgment

The financial support of the project (RISE-UAIC), code CNFIS-FDI-2026-F-0549, is acknowledged by I.D. Ghiba.

## Appendix A. The technical details of the proof of Theorem 2.12

In this section we give the proof of Theorem 2.12. Recall that $\hat{\boldsymbol{\tau}} = \hat{\boldsymbol{\tau}}(\log V)$. We prove only the implication from vector-monotonicity to matrix-monotonicity since the other implication is trivial. The first step in the proof is to show that for every $x \in \Pi$ and all $i \neq j$,

$$(\ell_i - \ell_j)(\hat{\tau}_i(\ell) - \hat{\tau}_j(\ell)) \geq 0. \tag{A.24}$$

In particular, $v$ and $g(\ell)$ are sorted in the same algebraic order.

Fix $\ell \in \Pi$ and indices $i \neq j$. Let $\pi$ denote the permutation exchanging only $i$ and $j$, and consider $h = \pi(\ell)$. Since permutations preserve the condition $\sum_{k=1}^n \ell_k = 0$, we still have $h \in \Pi$. By vector monotonicity we know that

$$0 \leq \langle \hat{\tau}(\ell) - \hat{\tau}(h),\, \ell - h \rangle. \tag{A.25}$$

By permutation-equivariance we have

$$\hat{\tau}_i(h) = \hat{\tau}_i(\ell_{\pi(1)}, \ell_{\pi(2)}, \dots, \ell_{\pi(n)}) = \hat{\tau}_{\pi(i)}(\ell_1, \ell_2, \dots, \ell_n). \tag{A.26}$$

Hence

$$0 \leq \langle \hat{\tau}(\ell) - \hat{\tau}(h),\, \ell - h \rangle = \sum_{i=1}^n (\hat{\tau}_i(\ell) - \hat{\tau}_{\pi(i)}(\ell))(\ell_i - \ell_{\pi(i)}). \tag{A.27}$$

Since the components which are not permuted lead to zero entries in the sum on the right-hand side of the above inequality, it follows that

$$0 \leq 2\,(\hat{\tau}_i(\ell) - \hat{\tau}_j(\ell))(\ell_i - \ell_j), \tag{A.28}$$

which proves the claim. As a consequence, if $\ell^{\downarrow}$ denotes $\ell$ sorted in decreasing order, then

$$\hat{\tau}(\ell)^{\downarrow} = \hat{\tau}(\ell^{\downarrow}). \tag{A.29}$$

The second step is the reduction of the matrix monotonicity to a diagonal–orthogonal form. Consider two arbitrary matrices in $\operatorname{dev}\operatorname{Sym}(n)$. Then, there exists $Q, R \in \mathrm{O}(n)$ and $\ell, h \in \Pi$ such that

$$X = Q\,\mathrm{diag}(\ell)\,Q^T, \qquad Y = R\,\mathrm{diag}(h)\,R^T. \tag{A.30}$$

Using orthogonal invariance of the trace, we obtain

$$\begin{aligned}\langle \hat{\boldsymbol{\tau}}(X) - \hat{\boldsymbol{\tau}}(Y),\, X - Y \rangle &= \langle \mathrm{diag}(\hat{\tau}(\ell)) - \widetilde{Q}\,\mathrm{diag}(\hat{\tau}(h))\,\widetilde{Q}^T,\, \mathrm{diag}(\ell) - \widetilde{Q}\,\mathrm{diag}(h)\,\widetilde{Q}^T \rangle \\ &= \langle \mathrm{diag}(\hat{\tau}(\ell)), \mathrm{diag}(\ell) \rangle + \langle \mathrm{diag}(\hat{\tau}(h)), \mathrm{diag}(h) \rangle \\ &\quad - \langle \mathrm{diag}(\hat{\tau}(\ell))\,\widetilde{Q}, \widetilde{Q}\,\mathrm{diag}(h) \rangle - \langle \mathrm{diag}(\ell)\,\widetilde{Q}, \widetilde{Q}\,\mathrm{diag}(\hat{\tau}(h)) \rangle,\end{aligned} \tag{A.31}$$

where $\widetilde{Q} := Q^T R \in O(3)$. We remark that the matrix $M = (m_{ij})$ given by $m_{ij} = \widetilde{Q}_{ij}^2$ is doubly stochastic, i.e., $m_{ij} \geq 0$, $\sum_j^n m_{ij} = 1$, $\sum_i^n m_{ij} = 1$, since the orthogonality of $\widetilde{Q}$ is exactly what guarantees that the row vectors and the column vectors have norm 1, and squaring the entries converts those norm-1 conditions into the row and column sum conditions for a doubly stochastic matrix. By the Birkhoff–von Neumann theorem (Birkhoff, 1946; Hardy et al., 1952; Horn, 1954), every doubly stochastic matrix can be written as a convex combination of permutation matrices. Hence there exist permutation matrices $\mathcal{P}_k$ and coefficients $q_k \geq 0$, $\sum_k q_k = 1$, such that $M = \sum_k q_k \mathcal{P}_k$. Therefore,

$$\langle a, M.b \rangle = \left\langle a, \left(\sum_k q_k \mathcal{P}_k\right).b \right\rangle = \sum_k q_k \langle a, \mathcal{P}_k.b \rangle. \tag{A.32}$$

Thus $\langle a, M.b \rangle$ is a convex combination of the numbers $\langle a, \mathcal{P}_k b \rangle$. Hence

$$\langle a, M.b \rangle \leq \max_{\mathcal{P} \text{ permutation matrix}} \langle a, \mathcal{P}.b \rangle. \tag{A.33}$$

If $P$ corresponds to a permutation $\pi$, then $\langle a, \mathcal{P}.b \rangle = \sum_i a_i b_{\pi(i)}$. We recall the rearrangement inequality (Hardy et al., 1952; Marshall et al., 2011): if $a_1 \leq a_2 \leq \cdots \leq a_n$, $b_1 \leq b_2 \leq \cdots \leq b_n$, then for every permutation $\pi$ of $\{1, 2, \dots, n\}$,

$$\sum_{i=1}^n a_i\, b_{3-i+1} \;\leq\; \sum_{i=1}^n a_i\, b_{\pi(i)} \;\leq\; \sum_{i=1}^n a_i\, b_i. \tag{A.34}$$

Therefore,

$$\max_{\mathcal{P} \text{ perm}} \langle a, \mathcal{P}.b \rangle = \langle a^{\downarrow}, b^{\downarrow} \rangle \qquad \forall\, a, b \in \mathbb{R}^n. \tag{A.35}$$

Combining the above steps, for all $\widetilde{Q} \in O(n)$, it holds that

$$\langle \mathrm{diag}(\ell)\widetilde{Q}, \widetilde{Q}.\,\mathrm{diag}(h)\rangle = \langle \mathrm{diag}(\ell), M.\,\mathrm{diag}(h)\rangle \le \max_{\mathcal{P}\ \mathrm{perm}} \langle \mathrm{diag}(\ell), \mathcal{P}.\,\mathrm{diag}(h)\rangle = \langle \mathrm{diag}(\ell)^{\downarrow}, \mathrm{diag}(h)^{\downarrow}\rangle, \tag{A.36}$$

$$\begin{aligned}\langle \mathrm{diag}(\hat\tau(\ell))\widetilde{Q}, \widetilde{Q}.\,\mathrm{diag}(\hat\tau(h))\rangle &= \langle \mathrm{diag}(\hat\tau(\ell)), M.\,\mathrm{diag}(\hat\tau(h))\rangle\\ &\le \max_{\mathcal{P}\ \mathrm{perm}} \langle \mathrm{diag}(\hat\tau(\ell)), \mathcal{P}.\,\mathrm{diag}(\hat\tau(h))\rangle\\ &= \langle \mathrm{diag}(\hat\tau(\ell))^{\downarrow}, \mathrm{diag}(\hat\tau(h))^{\downarrow}\rangle.\end{aligned} \tag{A.37}$$

Since inner products are invariant under simultaneous permutations, we have

$$\langle \mathrm{diag}(\hat\tau(\ell)), \mathrm{diag}(\ell)\rangle = \langle \hat\tau(\ell), \ell\rangle = \langle \hat\tau(\ell^{\downarrow}), \ell^{\downarrow}\rangle, \tag{A.38}$$

and similarly $\langle \mathrm{diag}(\hat\tau(h)), \mathrm{diag}(h)\rangle = \langle \hat\tau(h), h\rangle = \langle \hat\tau(h^{\downarrow}), h^{\downarrow}\rangle$.

We therefore obtain

$$\begin{aligned}\langle \hat{\boldsymbol\tau}(X) - \hat{\boldsymbol\tau}(Y),\ X - Y\rangle &\ge \langle \hat\tau(\ell^{\downarrow}), \ell^{\downarrow}\rangle + \langle \hat\tau(h^{\downarrow}), h^{\downarrow}\rangle\\ &\quad - \langle \hat\tau(\ell^{\downarrow}), h^{\downarrow}\rangle - \langle \ell^{\downarrow}, \hat\tau(h^{\downarrow})\rangle\\ &= \langle \hat\tau(\ell^{\downarrow}) - \hat\tau(h^{\downarrow}),\ \ell^{\downarrow} - h^{\downarrow}\rangle.\end{aligned} \tag{A.39}$$

Since $\ell^{\downarrow}, h^{\downarrow} \in \Pi$, vector monotonicity of $\hat\tau$ implies $\langle \hat\tau(\ell^{\downarrow}) - \hat\tau(h^{\downarrow}),\ \ell^{\downarrow} - h^{\downarrow}\rangle \ge 0$, which proves matrix monotonicity of $\hat{\boldsymbol\tau}$ on $\mathrm{dev}\,\mathrm{Sym}(n)$.

## Appendix B. An alternative proof for twice-differentiable energies

**Definition B.1** (*Second Derivative on $\Pi$*)**.** Let $\hat g^{\mathrm{inc}} : \Pi \to \mathbb{R}$ be twice-differentiable. The second derivative of $\hat g^{\mathrm{inc}}$ at $\ell \in \Pi$ is the mapping $\mathrm{D}^2\hat g^{\mathrm{inc}}(\ell) : \Pi \times \Pi \to \mathbb{R}$, defined by

$$\mathrm{D}^2\hat g^{\mathrm{inc}}(\ell).[h,k] := \frac{d}{dt}\Big|_{t=0} \mathrm{D}\hat g^{\mathrm{inc}}(\ell + th).\,k, \qquad h,k \in \Pi. \tag{B.1}$$

Equivalently, $\mathrm{D}^2\hat g^{\mathrm{inc}}(\ell)$ is a bilinear form on $\Pi$.

Let $\hat\tau$ be defined by $\langle \hat\tau(\ell), k\rangle = \mathrm{D}\hat g^{\mathrm{inc}}(\ell).\,k\ \forall \ell, k \in \Pi$. Fix $\ell \in \Pi$ and $h, k \in \Pi$, and consider the function $\phi(t) := \langle \hat\tau(\ell + th), k\rangle$. Using the defining relation, we have $\phi(t) = \mathrm{D}\hat g^{\mathrm{inc}}(\ell + th).\,k$. Differentiating with respect to $t$ at $t = 0$ yields $\phi'(0) = \langle \mathrm{D}\hat\tau(\ell).\,h, k\rangle = \mathrm{D}^2\hat g^{\mathrm{inc}}(\ell)[h,k]$. In particular,

$$\langle \mathrm{D}\hat\tau(\ell).\,h, h\rangle = \mathrm{D}^2\hat g^{\mathrm{inc}}(\ell).[h,h]. \tag{B.2}$$

**Remark B.2** (*Characterization by the Reduced Hessian on $\Pi$, see Remark 5.12 from Wollner et al. (2026a), too*)**.** Let $\hat g^{\mathrm{inc}} : \Pi \to \mathbb{R}$ be twice-differentiable. Define the linear map $A : \mathbb{R}^2 \to \Pi$, $A(x_1, x_2) := (x_1, x_2, -x_1 - x_2)$. In canonical coordinates, $A$ is represented by the matrix

$$A = \begin{pmatrix} 1 & 0\\ 0 & 1\\ -1 & -1\end{pmatrix}. \tag{B.3}$$

Define

$$\hat g^{\mathrm{inc}}_{\mathrm{red}} : \mathbb{R}^2 \to \mathbb{R}, \qquad \hat g^{\mathrm{inc}}_{\mathrm{red}}(x) := \hat g^{\mathrm{inc}}(A\,x) = \hat g^{\mathrm{inc}}(x_1, x_2, -x_1 - x_2). \tag{B.4}$$

If $\hat g^{\mathrm{inc}}$ admits a twice-differentiable extension to an open neighborhood in $\mathbb{R}^3$, denoted by $g$, then

$$\mathrm{D}^2\hat g^{\mathrm{inc}}_{\mathrm{red}}(x) = A^T \mathrm{D}^2 g(\ell)\,A, \qquad \ell = A\,x, \tag{B.5}$$

and the matrix $A^T\mathrm{D}^2 g(\ell)\,A$ is precisely the reduced Hessian of $g$ on the constraint manifold $\Pi$.

**Proof.** Since $A$ is linear and maps $\mathbb{R}^2$ bijectively onto $\Pi$, every $\ell \in \Pi$ can be written uniquely as $\ell = A\,x$, and every $h \in \Pi$ uniquely as $h = A\,\eta$. Now let

$$\hat g^{\mathrm{inc}}_{\mathrm{red}}(x) = \hat g^{\mathrm{inc}}(A\,x). \tag{B.6}$$

By the chain rule, for every $\eta \in \mathbb{R}^2$ one has

$$\mathrm{D}\hat g^{\mathrm{inc}}_{\mathrm{red}}(x).\eta = \mathrm{D}\hat g^{\mathrm{inc}}(A\,x).[A\,\eta]. \tag{B.7}$$

Since $A$ is linear, differentiating once more in the direction $\eta$ gives

$$\mathrm{D}^2\hat g^{\mathrm{inc}}_{\mathrm{red}}(x).[\eta,\eta] = \mathrm{D}^2\hat g^{\mathrm{inc}}(A\,x).[A\,\eta, A\,\eta]. \tag{B.8}$$

Thus, if $\ell = A\,x$ and $h = A\,\eta$, then

$$\mathrm{D}^2\hat g^{\mathrm{inc}}(\ell).[h,h] = \mathrm{D}^2\hat g^{\mathrm{inc}}_{\mathrm{red}}(x).[\eta,\eta] = \langle \eta, \mathrm{D}^2\hat g^{\mathrm{inc}}_{\mathrm{red}}(x).\eta\rangle. \tag{B.9}$$

Since $A$ is an isomorphism between $\mathbb{R}^2$ and $\Pi$, the condition $h \ne 0$ is equivalent to $\eta \ne 0$. Therefore,

$$\mathrm{D}^2\hat g^{\mathrm{inc}}(\ell).[h,h] > 0 \quad \forall\, h \in \Pi \setminus \{0\} \tag{B.10}$$

if and only if

$$\mathrm{D}^2\hat g^{\mathrm{inc}}_{\mathrm{red}}(x).[\eta,\eta] > 0 \quad \forall\, \eta \in \mathbb{R}^2 \setminus \{0\}, \tag{B.11}$$

which is equivalent to the positive definiteness of $\mathrm{D}^2\hat g^{\mathrm{inc}}_{\mathrm{red}}(x)$.

Finally, if $\hat g^{\mathrm{inc}}$ is extended to $g$ defined on an open neighborhood in $\mathbb{R}^3$, then the usual chain rule for Hessians yields

$$\mathrm{D}^2\hat g^{\mathrm{inc}}_{\mathrm{red}}(x) = A^T\mathrm{D}^2 g(\ell)\,A, \qquad \ell = A\,x. \tag{B.12}$$

This is exactly the Hessian of $g$ restricted to admissible directions $h \in \Pi$, hence it is the reduced Hessian associated with the constraint $\ell_1 + \ell_2 + \ell_3 = 0$. ■

**Lemma B.3** (*Positive semidefinite reduced Hessian implies monotonicity on $\Pi$*)**.** *Let $\hat g^{\mathrm{inc}} : \Pi \to \mathbb{R}$ be twice-differentiable. Define its intrinsic gradient $\hat\tau(\ell) \in \Pi$ by*

$$\langle \hat\tau(\ell), h\rangle = \mathrm{D}\hat g^{\mathrm{inc}}(\ell).\,h \qquad \forall\, \ell, h \in \Pi. \tag{B.13}$$

*Assume that the reduced Hessian is positive semidefinite. Then $\hat\tau$ is monotone on $\Pi$.*

**Proof.** Let $\ell, w \in \Pi$ be arbitrary. Since $\Pi$ is a linear subspace, the segment $\gamma(t) := w + t(\ell - w)$, $t \in [0,1]$, lies entirely in $\Pi$. Define

$$\varphi(t) := \langle \hat\tau(\gamma(t)),\ \ell - w\rangle. \tag{B.14}$$

Since $\hat g^{\mathrm{inc}}$ is twice-differentiable, the map $\hat\tau$ is of class differentiable on $\Pi$, and thus $\varphi$ is differentiable. By the chain rule,

$$\varphi'(t) = \left\langle \mathrm{D}\hat\tau(\gamma(t)).[\ell - w],\ \ell - w\right\rangle. \tag{B.15}$$

Because $\hat\tau$ is the intrinsic gradient of $\hat g^{\mathrm{inc}}$, its derivative is the intrinsic Hessian of $\hat g^{\mathrm{inc}}$, hence

$$\left\langle \mathrm{D}\hat\tau(\gamma(t)).\,h,\ k\right\rangle = \mathrm{D}^2\hat g^{\mathrm{inc}}(\gamma(t)).[h,k] \qquad \forall\, h,k \in \Pi. \tag{B.16}$$

In particular,

$$\varphi'(t) = \mathrm{D}^2\hat g^{\mathrm{inc}}(\gamma(t)).[\ell - w, \ell - w]. \tag{B.17}$$

By the assumed positive semidefiniteness of the reduced Hessian, we have $\varphi'(t) \ge 0\ \forall\, t \in [0,1]$. Therefore $\varphi$ is nondecreasing on $[0,1]$, and thus $\varphi(1) - \varphi(0) \ge 0$. This gives $\langle \hat\tau(\ell) - \hat\tau(w),\ \ell - w\rangle \ge 0$. Hence $\hat\tau$ is monotone on $\Pi$. ■

**Proposition B.4** (*Polyconvexity Implies Hill's Inequality in $\mathrm{SL}(2)$ for Objective Isotropic Energies*)**.** *Let $W^{\mathrm{inc}} : \mathrm{SL}(2) \to \mathbb{R}$ be an objective and isotropic elastic twice-differentiable energy. Assume that $W^{\mathrm{inc}}$ is polyconvex. Then $W^{\mathrm{inc}}$ satisfies the weak Hill inequality.*

**Proof.** Since $W^{\mathrm{inc}}$ is objective and isotropic, the energy depends only on the singular values of $F$. On $\mathrm{SL}(2)$ these are of the form $(\lambda, \lambda^{-1})$ with $\lambda \ge 1$. Hence there exists a function $\phi : [0,\infty) \to \mathbb{R}$ such that $W^{\mathrm{inc}}(F) = \phi\left(\lambda_{\max} - \frac{1}{\lambda_{\max}}\right)$. By Theorem 2.3, polyconvexity of $W^{\mathrm{inc}}$ is equivalent to the fact that $\phi$ is convex and nondecreasing on $[0,\infty)$.

If $X \in \operatorname{dev}\operatorname{Sym}(2)$, then $X$ admits a spectral decomposition

$$X = Q\,\operatorname{diag}(\ell)\,Q^T, \qquad Q \in \mathrm{O}(2),$$
$$\ell \in \Pi := \{v = (\ell_1, \ell_2) \in \mathbb{R}^2 \mid \ell_1 + \ell_2 = 0\}. \tag{B.18}$$

Define

$$\widehat{W}^{\mathrm{inc}}(X) := W^{\mathrm{inc}}(\exp(X)), \tag{B.19}$$

and let $\widehat{g}^{\mathrm{inc}} : \Pi \to \mathbb{R}$ be the corresponding reduced spectral function:

$$\widehat{W}^{\mathrm{inc}}(X) = \widehat{g}^{\mathrm{inc}}(\ell) \qquad \text{whenever } X = Q\,\operatorname{diag}(\ell)\,Q^T. \tag{B.20}$$

In dimension two, every $\ell \in \Pi$ has the form $\ell = (t, -t)$, $t \in \mathbb{R}$. By isotropy and symmetry, it is enough to consider $t \geq 0$. Then $\exp(X) = Q\,\operatorname{diag}(e^t, e^{-t})\,Q^T$, and therefore

$$\widehat{g}^{\mathrm{inc}}(t, -t) = \psi(t) := \phi(2\sinh t), \qquad t \geq 0. \tag{B.21}$$

We compute

$$\psi'(t) = 2\cosh t\,\phi'(2\sinh t), \tag{B.22}$$

and

$$\psi''(t) = 4\cosh^2 t\,\phi''(2\sinh t) + 2\sinh t\,\phi'(2\sinh t). \tag{B.23}$$

Since $\phi$ is convex and nondecreasing, we have $\phi'' \geq 0,\ \phi' \geq 0$ (in the classical sense when $\phi \in C^2$, or in the weak sense if $\phi$ is only convex). Since moreover $\cosh t > 0, \sinh t \geq 0$ for $t \geq 0$, it follows that

$$\psi''(t) \geq 0 \qquad \forall\, t \geq 0. \tag{B.24}$$

We now interpret this as positivity of the reduced Hessian of $\widehat{g}^{\mathrm{inc}}$ on $\Pi$. Let

$$\ell = (t, -t) \in \Pi, \qquad h = (s, -s) \in \Pi. \tag{B.25}$$

Since $\widehat{g}^{\mathrm{inc}}(\ell) = \psi(t)$, the second derivative of $\widehat{g}^{\mathrm{inc}}$ in the admissible direction $h$ is given by

$$\mathrm{D}^2\widehat{g}^{\mathrm{inc}}(\ell).[h, h] = \psi''(t)\,s^2 \geq 0. \tag{B.26}$$

Thus the reduced Hessian of $\widehat{g}^{\mathrm{inc}}$ is positive semidefinite on $\Pi$.

Let $\widehat{\tau}(\ell) \in \Pi$ be the intrinsic gradient of $\widehat{g}^{\mathrm{inc}}$ on $\Pi$. By the previous step, the reduced Hessian of $\widehat{g}^{\mathrm{inc}}$ is positive semidefinite. Hence, by the preceding Lemma B.3, the principal Kirchhoff stress mapping $\widehat{\tau} : \Pi \to \Pi$ is monotone and the rest of the proof follows from Theorem 2.12. ■

## Data availability

No data was used for the research described in the article.